\documentclass[11pt]{article}
\usepackage{amsmath}
\usepackage{amssymb}
\usepackage{latexsym}
\usepackage{graphicx}
\newcommand{\ds}{\displaystyle}
\newtheorem{theorem}{\bf Theorem}
\newtheorem{lemma}{\bf Lemma}

\begin{document}

\baselineskip 7mm

\vspace*{10mm}

\begin{center}
{\LARGE\bf 
Tail Probability and Singularity of Laplace-Stieltjes Transform 
of a Heavy Tailed Random Variable}\\

\vspace*{10mm}

\baselineskip 5mm

Kenji Nakagawa \\
Department of Electrical Engineering, \\
Nagaoka University of Technology, \\
Nagaoka, Niigata 940-2188, Japan \\
E-mail nakagawa@nagaokaut.ac.jp \\
\today
\end{center}

\vspace*{5mm}

\begin{abstract}
In this paper, we will give a sufficient condition for a 
non-negative random variable $X$ to be heavy tailed by 
investigating the Laplace-Stieltjes transform of the 
probability distribution function. We focus on the relation 
between the singularity at the real point of the axis 
of convergence and the asymptotic decay of the tail 
probability. Our theorem is a kind of Tauberian theorems.
\end{abstract}



\vspace*{5mm}

\section{Introduction}
We consider the asymptotic decay of the {\it tail probability} 
$P(X>x)$ of a heavy tailed random variable $X$. A random variable 
$X$ is said to be {\it heavy tailed} if 
\begin{eqnarray}
\eta\equiv\lim_{x\to\infty}\frac{\log P(X>x)}{\log x}<0.\label{eqn:1}
\end{eqnarray}
$\eta$ is called the decay rate of $P(X>x)$.

Let $F(x)$ be the probability distribution function (pdf) 
of $X$, i.e., $F(x)=P(X\leq x)$. For example, if $F(x)=1-1/x,\ x\geq 1$, 
then $P(X>x)=1/x$, hence $X$ is heavy tailed with decay 
rate $\eta=-1$.

In this paper, we will give a sufficient condition for $X$ 
to be heavy tailed by analytic properties of the Laplace-Stieltjes (LS) 
transform of $F(x)$. The LS transform of $F(x)$ is defined by 
\begin{eqnarray}
\varphi(s)=\int_0^\infty e^{-sx}dF(x).\label{eqn:2}
\end{eqnarray}

In general, for a function $R(x)$, which is of bounded variation 
in the interval $0\leq x\leq c$ for any positive $c$, the LS transform 
$\Psi(s)=\int_0^\infty e^{-sx}dR(x),\ s=\sigma+i\tau$, is defined. 
If $\Psi(s)$ converges for $\sigma>\sigma_0$ and diverges for 
$\sigma<\sigma_0$, then $\sigma_0$ is said to be the {\it abscissa 
of convergence} of $\Psi(s)$. The line $\Re s=\sigma_0$ is called 
the {\it axis of convergence}. In the 
case $\sigma_0=0$, some local information of $\Psi(s)$ 
at $s=0$ provides the asymptotic behavior of $R(x)$ as 
$x\to\infty$. Such a proposition is called a Tauberian 
theorem. The following is one of the Tauberian theorems.

\noindent{\bf Theorem} (Widder\,\cite{wid}, p.192, Theorem 4.3) 
{\it Let $R(x)$ be a non-decreasing function and the abscissa of 
convergence of $\Psi(s)$ be $\sigma_0=0$. If for constants $r\leq0$ 
and $A$
\begin{eqnarray}
\lim_{s\to0+}|\Psi(s)s^{-r}-A|=0,\label{eqn:3}
\end{eqnarray}
then
\begin{eqnarray}
\lim_{x\to\infty}|R(x)\Gamma(r+1)x^r-A|=0,\label{eqn:4}
\end{eqnarray}
where $\Gamma$ denotes the gamma function.}

In \cite{nak1},\cite{nak3},\cite{nak4}, we studied the 
asymptotic decay of a light tailed random variable. 
A random variable $X$ is said to be {\it light tailed} 
if the tail probability $P(X>x)$ decays exponentially, 
i.e., 
\begin{eqnarray}
\lim_{x\to\infty}\frac{1}{x}\log P(X>x)<0.\label{eqn:5}
\end{eqnarray}
We obtained in \cite{nak4} the following theorem which gives 
a sufficient condition for a light tailed random variable. 

\noindent{\bf Theorem} (Nakagawa\,\cite{nak4}) {\it For a 
non-negative random variable $X$ with probability distribution 
function $F(x)$, let $\varphi(s)=\int_0^\infty e^{-sx}dF(x)$ be 
the Laplace-Stieltjes transform of $F(x)$ and $\sigma_0$ be 
the abscissa of convergence of $\varphi(s)$. We assume 
$-\infty<\sigma_0<0$. If $s=\sigma_0$ is a pole of $\varphi(s)$, 
then we have}
\begin{eqnarray}
\lim_{x\to\infty}\frac{1}{x}\log P(X>x)=\sigma_0.\label{eqn:6}
\end{eqnarray}

It is known that for a monotonic $R(x)$, $s=\sigma_0$ is a 
singularity of $\Psi(s)$ (see Widder\,\cite{wid}, p.58, Theorem 5b). 
Since pdf $F(x)$ is monotonic increasing, $\sigma_0$ is 
a singularity of $\varphi(s)$. If $F(x)$ is the pdf of a heavy 
tailed random variable, then the abscissa of convergence 
of $\varphi(s)$ is necessarily $\sigma_0=0$ (see Widder\,\cite{wid}, 
p.40, Theorem 2.2b). Since $\varphi(0)=\int_0^\infty dF(x)=1$, 
$\sigma_0$ is not a pole, but other type of singularity. 

In the research of the asymptotic decay of a tail probability, 
we would like to construct a general theory such that a 
local analytic information of $\varphi(s)$ at $s=\sigma_0$ 
tells the asymptotic evaluation of the tail probability. 
In a light tailed case \cite{nak1},\cite{nak3},\cite{nak4}, we 
applied to this problem Ikehara's Tauberian theorem 
\cite{ike},\cite{kor} and its extension Graham-Vaaler's Tauberian 
theorem \cite{gra},\cite{kor}. Ikehara's theorem 
assumes a global analytic property of $\varphi(s)$, that is, 
$s=\sigma_0$ is a pole and there exist no other 
singularities on the axis of convergence $\Re s=\sigma_0$. 
While, Graham-Vaarler's theorem only assumes $s=\sigma_0$ 
is a pole, which yields weaker assertion than Ikehara's 
theorem, however, it is enough for our purpose to 
investigate the asymptotic decay of a tail probability. 
In this paper, we will apply Graham-Vaaler's Tauberian 
theorem to the decay of heavy tailed random variable. 

The purpose of our study is to apply our theory to the performance 
evaluation of the packet network. According to the 
research results for the Internet packet stream, it is reported 
that many characteristics are approximated by heavy tailed random 
variables \cite{cro}. Contrary to the light tailed case, if, for
example, the packet length is heavy tailed, the network performance, 
such as packet loss probability or end-to-end packet delay becomes 
worse. So, it is very important to investigate the tail probability 
of heavy tailed random variable from the view point of network engineering. 
When we apply the queueing theory to network engineering, 
even if we do not obtain the tail probability of packet 
length explicitly, we may obtain its Laplace transform by 
algebraic manipulation as Pollaczek-Khinchin formula \cite{kle}. Then, 
we know the singularity of the Laplace transform and the 
asymptotic decay can be investigated by our theory. 

Throughout this paper, we will use the following symbols. 
${\mathbb N}$, ${\mathbb N^+}$, ${\mathbb R}$, ${\mathbb C}$, 
denote the set of natural numbers, positive natural 
numbers, real numbers, complex numbers, and further, $\Re$, 
${\cal L}$, ${\cal F}$ denote the real part of a complex 
number, Laplace transform and Fourier transform, respectively.

\section{Examples of Heavy Tailed Random Variable}
Now, we look at some examples of heavy tailed random 
variables and their LS transforms. 
\subsection{Continuous Random Variable I}
Let $X$ be a random variable with pdf $F(x)=1-1/x,\ x\geq1$. 
The decay rate is $\eta=-1$. The LS transform $\varphi(s)$ 
of $F(x)$ is represented in a neighborhood of $s=0$ 
as 
\begin{align}
\varphi(s)&=\int_1^\infty e^{-sx}dF(x)\label{eqn:7}\\
&=s\log s+\beta(s),
\end{align}
where $\beta(s)$ is analytic in a neighborhood of $s=0$ 
(see Lemma \ref{lem:1}). 

\subsection{Continuous Random Variable II}
For a pdf $F(x)=1-1/\sqrt{x},\ x\geq 1$, we have the decay 
rate $\eta=-1/2$ and
\begin{eqnarray}
\varphi(s)=\ds\frac{2\pi}{\Gamma(1/2)}\sqrt{s}+\beta(s),\label{eqn:8}
\end{eqnarray}
where $\beta(s)$ is analytic in a neighborhood of 
$s=0$ (see Lemma \ref{lem:2}). 

\subsection{Discrete Random Variable}
Let $X$ be a discrete random variable with probability 
distribution $p=(p_n)_{n\in\mathbb N}$;
\begin{eqnarray}
p_n=\ds\frac{1}{(n+1)^r}-\ds\frac{1}{(n+2)^r},\ n\in\mathbb N,\ 
r\in\mathbb N^+.\label{eqn:9}
\end{eqnarray}
The tail probability of $X$ is $P(X>x)=\ds\frac{1}{(n+2)^r}$ 
and the decay rate is $\eta=-r$. The probability generating 
function (pgf) $p(z)$ of $p$ is
\begin{eqnarray}
p(z)=\ds\frac{z-1}{z^2}\sum_{n=1}^\infty\ds\frac{z^n}{n^r}
+\ds\frac{1}{z}.\label{eqn:10}
\end{eqnarray}
Substituting $z=e^{-s}$, we have
\begin{align}
\varphi(s)&=p(e^{-s})\label{eqn:11}\\
&=\alpha(s)s^r\log s+\beta(s),\label{eqn:12}
\end{align}
where $\alpha(s)$ and $\beta(s)$ are analytic in a 
neighborhood of $s=0$ with $\alpha(0)\neq0$ (see 
Lemma \ref{lem:3}).

\subsection{Stationary Distribution of M/G/1 Type Markov Chain}
We consider the tail probability of the stationary 
distribution of an M/G/1 type Markov chain \cite{fal}. 
Let $X$ be a heavy tailed discrete random variable with 
probability distribution $b=(b_n)_{n\in\mathbb N}$. 
Denote by $B(z)$ the pgf 
of $b$. Consider an M/G/1 type Markov chain with probability 
transition matrix
\begin{eqnarray}
P=
\left(\begin{array}{ccccc}
b_0 & b_1 & b_2 & b_3 & \ldots \\
a_0 & a_1 & a_2 & a_3 & \ldots \\
0   & a_0 & a_1 & a_2 & \ldots \label{eqn:13} \\
0   & 0   & a_0 & a_1 & \ldots \\
\vdots & \vdots & \vdots & \vdots & \ddots
\end{array}\right),\label{eqn:14}
\end{eqnarray}
where $a=(a_n)_{n\in\mathbb N}$ is a probability 
distribution with pgf $A(z)$. Suppose there exists the 
stationary distribution $\pi=(\pi_n)_{n\in\mathbb N}$ of $P$. 
Let $\pi(z)=\sum_{n=0}^\infty\pi_nz^n$ be the pgf 
of $\pi$, then by the Pollaczek-Khinchin formula \cite{fal},\cite{hay},\cite{nak2},
\begin{eqnarray}
\pi(z)=\ds\frac{\pi_0\left(zB(z)-A(z)\right)}{z-A(z)}.\label{eqn:15}
\end{eqnarray}
Substituting $z=e^{-s}$ into (\ref{eqn:15}), we have 
$\varphi(s)\equiv\pi(e^{-s})$. If the singularity of 
$\pi(z)$ at $z=1$ comes from the singularity of $B(z)$, 
we know the singularity of $\varphi(s)$ at $s=0$ comes from 
$B(e^{-s})$. So, 
it is expected that $\pi$ should be heavy tailed.

\section{Main Theorem}
From above examples, we expect that the following theorems 
hold. These are main theorems in this paper.
\begin{theorem}
\label{theo:1}
Let $X$ be a non-negative random variable with probability 
distribution function $F(x)$, and $\varphi(s)$ be the Laplace-Stieltjes 
transform of $F(x)$. Assume the abscissa of convergence of 
$\varphi(s)$ is $\sigma_0=0$, and $\varphi(s)$ is represented 
in a neighborhood of $s=0$ as 
\begin{eqnarray}
\varphi(s)=\alpha(s)s^r\log s+\beta(s),\ r\in{\mathbb N^+},\label{eqn:16}
\end{eqnarray}
where $\alpha(s)$ and $\beta(s)$ are analytic with 
$\alpha(0)\neq0$. Then, we have
\begin{eqnarray}
\lim_{x\to\infty}\ds\frac{\log P(X>x)}{\log x}=-r.\label{eqn:17}
\end{eqnarray}
\end{theorem}

\begin{theorem}
\label{theo:2}
Under the same notation as in Theorem $\ref{theo:1}$, if 
$\sigma_0=0$ and $\varphi(s)$ is represented in a neighborhood 
of $s=0$ as 
\begin{eqnarray}
\varphi(s)=\alpha(s)s^r+\beta(s),\ r>0,\ r\notin{\mathbb N},\label{eqn:18}
\end{eqnarray}
where $\alpha(s)$ and $\beta(s)$ are analytic with 
$\alpha(0)\neq0$. Then, we have
\begin{eqnarray}
\lim_{x\to\infty}\ds\frac{\log P(X>x)}{\log x}=-r.\label{eqn:19}
\end{eqnarray}
\end{theorem}

The following is a related work to Theorem 2.

\noindent{\bf Theorem} (Korevaar\,\cite{kor}, p.194, Theorem 8.2) 
{\it Let $S(x)$ vanish for $x<0$ and be locally integrable, positive 
and non-increasing for $x\geq0$. Let $\varphi(s)$ be the 
Laplace-Stieltjes transform of $S(x)$. Further, let $\rho(x)=x^{-r}l(x)$ 
with $0\leq r<1$ and slowly varying $l(x)$. Then, for a constant $A$, 
\begin{eqnarray}
\lim_{x\to\infty}\ds\frac{S(x)-S(\infty)}{\rho(x)}=A
\end{eqnarray}
if and only if
\begin{eqnarray}
\lim_{x\to\infty}\ds\frac{\varphi(1/x)-S(\infty)}{\rho(x)}=A\Gamma(1-r).
\end{eqnarray}
}

\section{Preliminary Lemmas}
We will prepare some lemmas for the proof of our main theorems. 
First, let us define the step function $\Delta_1(t)$ as
\begin{eqnarray}
\Delta_1(t)=\left\{\begin{array}{ll}\label{eqn:20}
1, & {\mbox{\rm if}}\ \ t\geq1,\\
0, & {\mbox{\rm if}}\ \ t<1.
\end{array}\right.
\end{eqnarray}

\begin{lemma}
\label{lem:1}
For $r\in{\mathbb N}^+$, we have
\begin{align}
\varphi(s)&\equiv{\cal L}\left(\ds\frac{1}{t^{r+1}}\Delta_1(t)\right)\label{eqn:21}\\
&=\int_1^\infty\ds\frac{1}{t^{r+1}}e^{-st}dt\label{eqn:22}\\
&=\ds\frac{(-1)^r}{r!}s^r\log s+\beta(s),\label{eqn:23}
\end{align}
where $\beta(s)$ is analytic in a neighborhood of $s=0$. 
\end{lemma}

\noindent{\bf Proof}\ \ The $(r+1)$th derivative of 
$\varphi(s)$ is 
\begin{align}
\varphi(s)^{(r+1)}(s)&=(-1)^{r+1}\int_1^\infty e^{-st}dt\label{eqn:24}\\
&=(-1)^{r+1}\left(\ds\frac{1}{s}+\ds\frac{e^{-s}-1}{s}\right).\label{eqn:25}
\end{align}
Since $(e^{-s}-1)/s$ is analytic, we have, by successive 
integrations, the desired result.\hfill$\Box$

\begin{lemma}
\label{lem:2}
For $r>0,\ r\notin{\mathbb N}$, let $r_0=\lfloor r\rfloor$ be 
the maximum integer among integers smaller than $r$, and 
let $\bar{r}=r-r_0$. Then we have 
\begin{align}
\varphi(s)&\equiv{\cal L}\left(\ds\frac{1}{t^{r+1}}\Delta_1(t)\right)\label{eqn:26}\\
&=\ds\frac{(-1)^{r_0+1}\pi}{\Gamma(r+1)\sin\pi\bar{r}}s^r+\beta(s),\label{eqn:27}
\end{align}
where $\beta(s)$ is analytic in a neighborhood of $s=0$.
\end{lemma}

\noindent{\bf Proof}\ \ By a formula of Laplace transform 
(see \cite{mor}, p.287),
\begin{align}
\varphi^{(r_0+1)}(s)&=(-1)^{r_0+1}\int_1^\infty t^{-\bar{r}}e^{-st}dt\label{eqn:28}\\
&=(-1)^{r_0+1}\left\{\int_0^\infty-\int_0^1\right\}t^{-\bar{r}}e^{-st}dt\label{eqn:29}\\
&=(-1)^{r_0+1}\Gamma(1-\bar{r})s^{\bar{r}-1}+\tilde\beta(s),\label{eqn:30}
\end{align}
where $\tilde\beta(s)$ is analytic in a neighborhood of $s=0$. 
By successive integrations,
\begin{align}
\varphi(s)&=(-1)^{r_0+1}\Gamma(1-\bar{r})\ds\frac{1}{\bar{r}}\cdot
\ds\frac{1}{\bar{r}+1}\cdot\ldots\cdot\ds\frac{1}{\bar{r}+r_0}s^{\bar{r}+r_0}+\beta(s)\label{eqn:31}\\
&=\ds\frac{(-1)^{r_0+1}\pi}{\Gamma(r+1)\sin\pi\bar{r}}s^r+\beta(s).\label{eqn:32}
\end{align}
In (\ref{eqn:32}), we applied the formulas; 
$\Gamma(z+1)=z\Gamma(z)$ and $\Gamma(z)\Gamma(1-z)=\pi/\sin\pi z$.\hfill$\Box$

\begin{lemma}
\label{lem:3}
Let $X$ be a discrete random variable with probability 
distribution $p=(p_n)_{n\in{\mathbb N}};$
\begin{eqnarray}
p_n=\ds\frac{1}{(n+1)^r}-\ds\frac{1}{(n+2)^r},\ n\in{\mathbb N},\ r\in{\mathbb N^+},\label{eqn:33}
\end{eqnarray}
and let $p(z)=\sum_{n=0}^\infty p_nz^n$ be the pgf of $p$. Then, $\varphi(s)=p(e^{-s})$ is 
represented in a neighborhood of $s=0$ as 
\begin{eqnarray}
\varphi(s)=\alpha(s)s^r\log s+\beta(s),\label{eqn:34}
\end{eqnarray}
where $\alpha(s)$ and $\beta(s)$ are analytic 
with $\alpha(0)\neq0$.
\end{lemma}

\noindent{\bf Proof}\ \ For $r\geq2$, by calculation 
and Riemann's formula (see Widder\,\cite{wid}, p.232), 
\begin{align}
p(z)&=\ds\frac{z-1}{z^2}\sum_{n=1}^\infty\frac{z^n}{n^r}+\frac{1}{z}\label{eqn:35}\\
&=\ds\frac{1}{(r-1)!}\frac{z-1}{z^2}\int_0^\infty\frac{zt^{r-1}}{e^t-z}dt
+\frac{1}{z}.\label{eqn:36}
\end{align}
Then, 
\begin{eqnarray}
\varphi(s)=\ds\frac{e^s(1-e^s)}{(r-1)!}\int_0^\infty\frac{t^{r-1}}{e^{t+s}-1}dt+e^s.\label{eqn:36+1}
\end{eqnarray}
By the change of variable $u=e^{t+s}-1$, we have for sufficiently 
small $s>0$, 
\begin{align}
\varphi(s)&=\ds\frac{e^s(1-e^s)}{(r-1)!}\int_{e^s-1}^\infty\frac{\{\log(1+u)-s\}^{r-1}}{u(u+1)}du
+\mbox{\rm regular\ term}\label{eqn:37}\\
&=\ds\frac{e^s(1-e^s)}{(r-1)!}\int_{e^s-1}\left(\ds\frac{1}{u}-\frac{1}{u+1}\right)
\left\{(-s)^{r-1}+(r-1)(-s)^{r-2}u+\ldots\right\}du+\mbox{\rm r.t.}\label{eqn:38}\\
&=\ds\frac{e^s(1-e^s)}{(r-1)!}(-s)^{r-1}\int_{e^s-1}\ds\frac{1}{u}+\mbox{\rm r.t.}\label{eqn:39}\\
&=\alpha(s)s^r\log s+\beta(s),\label{eqn:40}
\end{align}
where
\begin{eqnarray}
\alpha(s)=\ds\frac{(-1)^{r+1}}{(r-1)!}e^s\frac{e^s-1}{s},\ 
\alpha(0)=\frac{(-1)^{r+1}}{(r-1)!}\neq0\label{eqn:41}
\end{eqnarray}
and $\beta(s)$ is some analytic function.

For $r=1$, we have
\begin{eqnarray}
p(z)=\ds\frac{1-z}{z^2}\log(1-z)+\frac{1}{z}.\label{eqn:42}
\end{eqnarray}
Similar argument leads to the desired result. \hfill$\Box$

\begin{lemma}
\label{lem:4}
For a pdf $F(x)$, let 
\begin{eqnarray}
\varphi(s)=\int_0^\infty e^{-sx}dF(x)\label{eqn:43}
\end{eqnarray}
have the abscissa of convergence $\sigma_0=0$. If, in a 
neighborhood of $s=0$, 
\begin{eqnarray}
\varphi(s)=\alpha(s)s^r\log s+\beta(s),\ r\in{\mathbb N}^+,\label{eqn:44}
\end{eqnarray}
where $\alpha(s)$, $\beta(s)$ are analytic with $\alpha(0)\neq0$, 
then, 
\begin{eqnarray}
\alpha(0)\left\{\begin{array}{ll}\label{eqn:45}
>0, & \mbox{\rm if}\ r\ \mbox{\rm is\ odd},\\
<0, & \mbox{\rm if}\ r\ \mbox{\rm is\ even}.
\end{array}\right.
\end{eqnarray}
\end{lemma}

\noindent{\bf Proof}\ \ From (\ref{eqn:43}), 
\begin{eqnarray}
\varphi^{(r)}(0+)\left\{\begin{array}{ll}\label{eqn:46}
\leq0, & \mbox{\rm if}\ r\ \mbox{\rm is\ odd},\\
\geq0, & \mbox{\rm if}\ r\ \mbox{\rm is\ even},
\end{array}\right.
\end{eqnarray}
while from (\ref{eqn:44}) by calculation
\begin{eqnarray}
\varphi^{(r)}(s)=r!\alpha(s)\log s+\theta(s),\label{eqn:47}
\end{eqnarray}
where $\theta(s)$ is a function of $s$ with $|\theta(0+)|<\infty$. 
Thus, by (\ref{eqn:47})
\begin{eqnarray}
\varphi^{(r)}(0+)=r!\alpha(0)\times(-\infty)+\theta(0+).\label{eqn:48}
\end{eqnarray}
Comparing (\ref{eqn:46}), (\ref{eqn:48}), we have the result. 
\hfill$\Box$

\begin{lemma}
\label{lem:5}
Under the same notation as in Lemmas $\ref{lem:2}$, $\ref{lem:4}$, if 
$\varphi(s)$ has the abscissa of convergence $\sigma_0=0$ and 
\begin{eqnarray}
\varphi(s)=\alpha(s)s^r+\beta(s),\ r>0,\ r\notin{\mathbb N},\label{eqn:49}
\end{eqnarray}
with $\alpha(0)\neq0$, then, 
\begin{eqnarray}
\alpha(0)\left\{\begin{array}{ll}\label{eqn:50}
>0, & \mbox{\rm if}\ r_0=\lfloor r\rfloor\ \mbox{\rm is\ odd},\\
<0, & \mbox{\rm if}\ r_0\ \mbox{\rm is\ even}.
\end{array}\right.
\end{eqnarray}
\end{lemma}

\noindent{\bf Proof}\ \ From (\ref{eqn:43})
\begin{eqnarray}
\varphi^{(r_0+1)}(0+)\left\{\begin{array}{ll}\label{eqn:51}
\geq0, & \mbox{\rm if}\ r_0\ \mbox{\rm is\ odd},\\
\leq0, & \mbox{\rm if}\ r_0\ \mbox{\rm is\ even},
\end{array}\right.
\end{eqnarray}
while from (\ref{eqn:49})
\begin{eqnarray}
\varphi^{(r_0+1)}(s)=\sum_{k=0}^{r_0+1}\binom{r_0+1}{k}\alpha^{(r_0+1-k)}(s)
\cdot\frac{d}{ds^k}s^r+\beta^{(r_0+1)}(s).\label{eqn:52}
\end{eqnarray}
Since
\begin{eqnarray}
\frac{d}{ds^k}s^r\Big|_{s=0+}=\left\{\begin{array}{ll}\label{eqn:53}
0, & \mbox{\rm for}\ k=0,1,\ldots,r_0\\
+\infty, & \mbox{\rm for}\ k=r_0+1,
\end{array}\right.
\end{eqnarray}
we have from (\ref{eqn:52})
\begin{eqnarray}
\varphi^{(r_0+1)}(0+)=\alpha(0)\times(+\infty)+\beta^{(r_0+1)}(0).\label{eqn:54}
\end{eqnarray}
Comparing (\ref{eqn:51}), (\ref{eqn:54}), we have the 
result.\hfill$\Box$

\subsection{First Several Terms of $\alpha(s)$ and $\beta(s)$}
We will need later, in the proof of main theorems, 
to make Laplace transforms which have the same first 
several terms as those in the Taylor expansion of $\alpha(s)$ 
and $\beta(s)$, respectively.

First, consider the following case
\begin{eqnarray}
\varphi(s)=\alpha(s)s^r\log s+\beta(s),\ r\in{\mathbb N}^+,\label{eqn:55}
\end{eqnarray}
with expansions $\alpha(s)=\sum_{n=0}^\infty\alpha_nz^n$, 
$\alpha(0)\neq0$, and $\beta(s)=\sum_{n=0}^\infty\beta_nz^n$.

We will make functions $g^\ast(t)$ and $h^\ast(t)$ 
such that their Laplace transforms $G^\ast(s)\equiv{\cal L}(g^\ast(t))$ and 
$H^\ast(s)\equiv{\cal L}(h^\ast(t))$ satisfy the 
following (i) and (ii) for $L\in{\mathbb N}^+$.
\begin{itemize}
\item[(i)] $G^\ast(s)=\alpha^\ast(s)s^r\log s+\tilde\beta(s)$,\ where 
$\alpha^\ast(s)=\ds\sum_{n=0}^{L-1}\alpha_ns^n$ and $\tilde\beta(s)$ 
is analytic in a neighborhood of $s=0$. Let 
$\tilde\beta(s)=\ds\sum_{n=0}^\infty\tilde\beta_ns^n$ be the 
expansion at $s=0$.
\item[(ii)] $H^\ast(s)=\ds\sum_{n=0}^{L-1}(\beta_n-\tilde\beta_n)s^n+$
higher order terms.
\end{itemize}

The above (i) and (ii) mean that $\alpha^\ast(s)$ is equal to 
the sum of the first $L$ terms of $\alpha(s)$, and the sum of the 
first $L$ terms of $H^\ast(s)+\tilde\beta(s)$ is equal to that of 
$\beta(s)$. 

\noindent\underline{Function $g^\ast(t)$}

Define
\begin{eqnarray}
&&g^\ast(t)=\sum_{k=0}^{L-1}\frac{g_k}{t^{r+k+1}}\Delta_1(t),\ t\in{\mathbb R},\label{eqn:56}\\
&&g_k=(-1)^{r+k+1}(r+k)!\alpha_k,\ k=0,1,\ldots,L-1.\label{eqn:57}
\end{eqnarray}
where $\Delta_1(t)$ was defined in (\ref{eqn:20}). 
By Lemma \ref{lem:1}, we see that $g^\ast(t)$ satisfies (i). The first 
coefficient $g_0$ is positive by Lemma \ref{lem:4}, i.e.,
\begin{eqnarray}
g_0=(-1)^{r+1}r!\alpha_0>0.\label{eqn:58}
\end{eqnarray}

\noindent\underline{Function $h^\ast(t)$}

Let $h_k(t)=ke^{-kt},\ t\geq0,\ k=1,2,\ldots,$ and $H_k(s)={\cal L}(h_k(t))$. 
We have $H_k(s)=k/(s+k),\ \Re s>-k$, and the expansion
\begin{eqnarray}
H_k(s)=\sum_{n=0}^\infty\left(-\frac{s}{k}\right)^n,\ |s|<k.\label{eqn:59}
\end{eqnarray}
The sum of the first $L$ terms of (\ref{eqn:59}) is 
represented as
\begin{eqnarray}
\left(1,-\ds\frac{1}{k},\ldots,\left(-\ds\frac{1}{k}\right)^{L-1}\right)
\left(1,s,\ldots,s^{L-1}\right)^T,\label{eqn:60}
\end{eqnarray}
where $^T$ denotes the transposition of vector. Let $V$ be 
the $L\times L$ matrix;
\begin{eqnarray}
V=
\left(\begin{array}{cccc}\label{eqn:61}
1 & -1 & \ldots & \left(-1\right)^{L-1} \\
1 & -\ds\frac{1}{2} & \ldots & \left(-\ds\frac{1}{2}\right)^{L-1} \\
\vdots & \vdots & & \vdots \\
1 & -\ds\frac{1}{L} & \ldots & \left(-\ds\frac{1}{L}\right)^{L-1}
\end{array}\right).
\end{eqnarray}
We will make a desired function by a linear combination 
of $h_k(t),\ k=1,2,\ldots,L$. Let $h(t)=\sum_{k=1}^Ld_kh_k(t),\ t\geq0$, 
and write $\mbox{{\boldmath$d$}}=(d_1,d_2,\ldots,d_L)$, 
$\mbox{{\boldmath$s$}}=(1,s,\ldots,s^{L-1})$. Further, write 
$\mbox{{\boldmath$\beta$}}=(\beta_0,\beta_1,\ldots,\beta_{L-1})$, 
$\mbox{{\boldmath$\tilde\beta$}}=(\tilde\beta_0,\tilde\beta_1,\ldots,\tilde\beta_{L-1})$. 
The sum of the first $L$ terms of $H(s)={\cal L}(h(t))$ is 
{\boldmath$d$}\,V{\boldmath$s$}$^T$, then we must solve the 
equation 
\begin{eqnarray}
\mbox{{\boldmath$d$}\,V{\boldmath$s$}$^T$=
$(${\boldmath$\beta$}$-${\boldmath$\tilde\beta$}$)${\boldmath$s$}$^T$}.\label{eqn:62}
\end{eqnarray}
Since $\det V\neq0$ (Vandermonde matrix), we have 
$\mbox{{\boldmath$d$}}=(\mbox{{\boldmath$\beta$}}-{\mbox{\boldmath$\tilde\beta$}})V^{-1}$. 
We write this solution as 
$\mbox{{\boldmath$d$}}=(d_1,d_2,\ldots,d_L)$, 
then
\begin{eqnarray}
h^\ast(t)=\sum_{k=1}^Ld_kh_k(t)\label{eqn:63}
\end{eqnarray}
is a desired function, i.e., $H^\ast(s)={\cal L}(h^\ast(t))$ 
satisfies (ii). 

Summarizing above,

\begin{lemma}
\label{lem:6}
Let $\varphi(s)$ be the LS transform of a pdf and 
the abscissa of convergence be $\sigma_0=0$. If
\begin{eqnarray}
\varphi(s)=\alpha(s)s^r\log s+\beta(s),\ r\in{\mathbb N}^+,\label{eqn:64}
\end{eqnarray}
where $\alpha(s)$, $\beta(s)$ are analytic in a neighborhood of 
$s=0$ with $\alpha(0)\neq0$, then $g^\ast(t)$ in $(\ref{eqn:56})$ 
and $h^\ast(t)$ in $(\ref{eqn:63})$ satisfy $\mbox{\rm(i)}$ and $\mbox{\rm(ii)}$.
\end{lemma}

Similarly, in the case
\begin{eqnarray}
\varphi(s)=\alpha(s)s^r+\beta(s),\ r>0,\ r\notin{\mathbb N},\label{eqn:65}
\end{eqnarray}
with $\alpha(s)=\sum_{n=0}^\infty\alpha_ns^n$, $\alpha_0\neq0$, 
$\beta(s)=\sum_{n=0}^\infty\beta_ns^n$, we will make 
functions $g^\ast(t)$, $h^\ast(t)$ such that their Laplace 
transforms $G^\ast(s)$, $H^\ast(s)$ satisfy the following 
(iii) and (iv) for $L\in{\mathbb N}^+$. 
\begin{itemize}
\item[(iii)] $G^\ast(s)=\alpha^\ast(s)s^r+\tilde\beta(s)$, where 
$\alpha^\ast(s)=\ds\sum_{n=0}^{L-1}\alpha_ns^n$ and $\tilde\beta(s)$ 
is analytic in a neighborhood of $s=0$. Let 
$\tilde\beta(s)=\ds\sum_{n=0}^\infty\tilde\beta_ns^n$ be the 
expansion at $s=0$.
\item[(iv)] $H^\ast(s)=\ds\sum_{n=0}^{L-1}(\beta_n-\tilde\beta_n)s^n+$
higher order terms.
\end{itemize}
In this case,
\begin{eqnarray}
&&g^\ast(t)=\sum_{k=0}^{L-1}\frac{g_k}{t^{r+k+1}}\Delta_1(t),\ t\in{\mathbb R},\label{eqn:66}\\
&&g_k=(-1)^{r_0+k+1}\frac{\sin\pi\bar{r}}{\pi}\Gamma(r+k+1)\alpha_k,\ k=0,1,\ldots,L-1.\label{eqn:67}
\end{eqnarray}
satisfies (iii). The first 
coefficient $g_0$ is positive by Lemma \ref{lem:5}, i.e.,
\begin{eqnarray}
g_0=(-1)^{r_0+1}\ds\frac{\sin\pi\bar{r}}{\pi}\Gamma(r+1)\alpha_0>0.\label{eqn:67+1}
\end{eqnarray}
The same $h^\ast(t)$ as (\ref{eqn:63}) satisfies (iv). Thus, we have

\begin{lemma}
\label{lem:7}
Let $\varphi(s)$ be the LS transform of a pdf and the abscissa of 
convergence be $\sigma_0=0$. If 
\begin{eqnarray}
\varphi(s)=\alpha(s)s^r+\beta(s),\ r>0,\ r\notin{\mathbb N},\label{eqn:68}
\end{eqnarray}
then $g^\ast(t)$ in $(\ref{eqn:66})$ and $h^\ast(t)$ in 
$(\ref{eqn:63})$ satisfy $\mbox{\rm(iii)}$ and $\mbox{\rm(iv)}$.
\end{lemma}

\subsection{Majorant and Minorant Functions}
For the evaluation of the tail probability $P(X>x)$ from 
above and below, we need to use majorant and minorant 
functions for an exponential function (see Korevaar\,\cite{kor}, 
p.132, Graham-Vaaler\,\cite{gra}). If two functions $f_1$, $f_2$ 
satisfy $f_1(t)\geq f_2(t),\ t\in{\mathbb R}$, then $f_1$ is 
said to be a {\it majorant} for $f_2$, and $f_2$ is a {\it minorant} 
for $f_1$.

For $\omega>0$, we will define a majorant $M_\omega^1(t)$ 
and a minorant $m_\omega^1(t)$ for
\begin{eqnarray}
E_\omega(t)\equiv\left\{\begin{array}{ll}\label{eqn:69}
e^{-\omega t}, & t\geq0\\
0, & t<0.
\end{array}\right.
\end{eqnarray}
Define (see Korevaar\,\cite{kor}, p.132)
\begin{align}
M_\omega^1(t)&=\left(\frac{\sin\pi t}{\pi}\right)^2Q_\omega(t),
\ t\in{\mathbb R},\label{eqn:70}\\
Q_\omega(t)&=\sum_{n=0}^\infty\frac{e^{-n\omega}}{(t-n)^2}
-\omega\sum_{n=1}^\infty e^{-n\omega}\left(\frac{1}{t-n}-\frac{1}{t}\right),\label{eqn:71}
\end{align}
and 
\begin{align}
m_\omega^1(t)&=M_\omega^1(t)-\left(\frac{\sin\pi t}{\pi t}\right)^2,\ t\in{\mathbb R}.\label{eqn:72}
\end{align}
Moreover, for $L\in{\mathbb N}^+$, define
\begin{eqnarray}
M_\omega^L(t)=\left(M_\omega^1(t)\right)^L\ \ \mbox{\rm and}\ \ 
m_\omega^L(t)=\left(m_\omega^1(t)\right)^L.\label{eqn:73}
\end{eqnarray}
For $\sigma>0$, $\delta>0$, write $\omega=2\pi\sigma/\delta$, then define
\begin{eqnarray}
M_{\sigma,\delta}^L(t)\equiv M_\omega^L\left(\frac{\delta t}{2\pi}\right)=
M_{2\pi\sigma/\delta}^L\left(\frac{\delta t}{2\pi}\right),\label{eqn:74}\\
m_{\sigma,\delta}^L(t)\equiv m_\omega^L\left(\frac{\delta t}{2\pi}\right)=
m_{2\pi\sigma/\delta}^L\left(\frac{\delta t}{2\pi}\right).\label{eqn:75}
\end{eqnarray}

\begin{lemma}
\label{lem:8}
$(\mbox{\rm Korevaar\,\cite{kor}})$\ \ 
For any $L\in{\mathbb N}^+$, $\sigma>0$, $\delta>0$, 
\begin{eqnarray}
M_{\sigma,\delta}^L,\ m_{\sigma,\delta}^L\in L^1({\mathbb R})\cap L^2({\mathbb R}).\label{eqn:76}
\end{eqnarray}
\end{lemma}

For $\lambda>0$, an entire function $f(z)$ of a complex variable $z=x+iy$ 
is of {\it exponential type} $\lambda$ if 
\begin{eqnarray}
|f(z)|\leq C\exp(\lambda|z|),\ z\in{\mathbb C},\ C>0.\label{eqn:77}
\end{eqnarray}
A real function $f(x)$ is of {\it type} $\lambda$ if $f(x)$ is 
the restriction to $\mathbb R$ of an entire function 
of exponential type $\lambda$.

\begin{lemma}
\label{lem:9}
$(\mbox{\rm Korevaar\,\cite{kor},\ Nakagawa\,\cite{nak4}})$\ \ 
$M_{\sigma,\delta}^L$ and $m_{\sigma,\delta}^L$ are of type 
$L\delta$.
\end{lemma}

\begin{lemma}
\label{lem:10}
$(\mbox{\rm Korevaar\,\cite{kor},\ Graham-Vaaler\,\cite{gra}})$\ \ 
For $L\in{\mathbb N}^+$, 
\begin{eqnarray}
E_{L\sigma}(t)\leq M_{\sigma,\delta}^L(t),\ t\in{\mathbb R},\label{eqn:79}
\end{eqnarray}
and for odd $L\in{\mathbb N}^+$, 
\begin{eqnarray}
m_{\sigma,\delta}^L(t)\leq E_{L\sigma}(t),\ t\in{\mathbb R}.\label{eqn:80}
\end{eqnarray}
\end{lemma}

\noindent{\bf Proof}\ \ If $L=1$, the result follows 
Korevaar\,\cite{kor}, p.129, Proposition 5.2. The odd power 
preserves the order of real numbers. \hfill$\Box$

From Lemma \ref{lem:8}, we can define the Fourier transforms 
$\hat{M}_{\sigma,\delta}^L={\cal F}(M_{\sigma,\delta}^L)$ and 
$\hat{m}_{\sigma,\delta}^L={\cal F}(m_{\sigma,\delta}^L)$, where 
the Fourier transform is defined as
\begin{eqnarray}
\hat{M}_{\sigma,\delta}^L(\tau)=\int_{-\infty}^\infty M_{\sigma,\delta}^L(t)
e^{-i\tau t}dt,\ \tau\in{\mathbb R}.\label{eqn:81}
\end{eqnarray}
Then, from Lemma \ref{lem:9} and the Paley-Wiener theorem, we have 

\begin{lemma}
\label{lem:11}
$(\mbox{\rm Rudin\,\cite{rud},\ Korevaar\,\cite{kor}})$\ \ For any 
$L\in{\mathbb N}^+$, 
\begin{eqnarray}
\mbox{\rm supp}(\hat{M}_{\sigma,\delta}^L)\subset[-L\delta, L\delta]\ \ \ 
\mbox{\rm and}\ \ \ \mbox{\rm supp}(\hat{m}_{\sigma,\delta}^L)
\subset[-L\delta, L\delta],\label{eqn:82}
\end{eqnarray}
where {\rm supp} denotes the support of a function. 
\end{lemma}

\subsection{Calculation of $\hat{M}_{\sigma,\delta}^L$ and 
$\hat{m}_{\sigma,\delta}^L$}
It is not difficult to calculate the Fourier transforms 
$\hat{M}_{\omega}^1={\cal F}(M_{\omega}^1)$ and 
$\hat{m}_{\omega}^1={\cal F}(m_{\omega}^1)$. 

Define 
\begin{eqnarray}
q_1(t)=\left(\frac{\sin\pi t}{\pi t}\right)^2,\ \ 
q_2(t)=\frac{\sin^2\pi t}{\pi t},\ t\in{\mathbb R},\label{eqn:83}
\end{eqnarray}
and write $\hat{q}_1={\cal F}(q_1)$, $\hat{q}_2={\cal F}(q_2)$. By 
calculation, we have
\begin{eqnarray}\label{eqn:84}
\hat{q}_1(\tau)=\left\{\begin{array}{ll}
1+\ds\frac{\tau}{2\pi}, & -2\pi\leq\tau<0,\\[4mm]
1-\ds\frac{\tau}{2\pi}, & 0\leq\tau<2\pi,\\[4mm]
0, & \mbox{\rm otherwise}.
\end{array}\right.
\end{eqnarray}
and
\begin{eqnarray}\label{eqn:85}
\hat{q}_2(\tau)=\left\{\begin{array}{ll}
\ds\frac{i}{2}, & -2\pi\leq\tau<0,\\[4mm]
-\ds\frac{i}{2}, & 0\leq\tau<2\pi,\\[4mm]
0, & \mbox{\rm otherwise}.
\end{array}\right.
\end{eqnarray}

\begin{lemma}
\label{lem:12}
We have
\begin{eqnarray}
\hat{M}_\omega^1(\tau)=\ds\frac{1}{1-e^{-(\omega+i\tau)}}\hat{q}_1(\tau)
-\frac{\omega}{\pi}\left(\frac{1}{1-e^{-(\omega+i\tau)}}-\frac{1}{1-e^{-\omega}}
\right)\hat{q}_2(\tau),\label{eqn:86}
\end{eqnarray}
and
\begin{eqnarray}
\hat{m}_\omega^1(\tau)=\ds\frac{e^{-(\omega+i\tau)}}{1-e^{-(\omega+i\tau)}}\hat{q}_1(\tau)
-\frac{\omega}{\pi}\left(\frac{1}{1-e^{-(\omega+i\tau)}}-\frac{1}{1-e^{-\omega}}
\right)\hat{q}_2(\tau).\label{eqn:87}
\end{eqnarray}
\end{lemma}

\noindent{\bf Proof}\ \ See Appendix \ref{app:A}.\hfill$\Box$

\medskip

Next, we will calculate $\hat{M}_\omega^L={\cal F}(M_\omega^L)$, 
$\hat{m}_\omega^L={\cal F}(m_\omega^L)$ and then calculate 
$\lim_{\omega\to0+}\hat{M}_\omega^L(\tau)$, $\lim_{\omega\to0+}\hat{m}_\omega^L(\tau)$, 
for $\tau\neq0$.

Let us define 
\begin{align}
u_\omega(t)&=\left(\frac{\sin\pi t}{\pi t}\right)^2
-\frac{\omega}{\pi}\frac{\sin^2\pi t}{\pi t}\label{eqn:88}\\[2mm]
&=q_1(t)-\frac{\omega}{\pi}q_2(t),\ t\in{\mathbb R},\label{eqn:89}\\[2mm]
v_\omega(t)&=\frac{\omega}{\pi}\frac{\sin^2\pi t}{\pi t}\label{eqn:90}\\[2mm]
&=\frac{\omega}{\pi}q_2(t),\ t\in{\mathbb R},\label{eqn:91}
\end{align}
and $\hat{u}_\omega={\cal F}(u_\omega)$, $\hat{v}_\omega={\cal F}(v_\omega)$.

\begin{lemma}
\label{lem:13}
We have
\begin{align}
\hat{M}_\omega^L(\tau)=\frac{1}{(2\pi)^{L-1}}\sum_{l=0}^L
\binom{L}{l}\left(\frac{1}{1-e^{-(\omega+i\tau)}}\right)^l\left(\frac{1}{1-e^{-\omega}}\right)^{L-l}
\hat{u}_\omega^{\ast l}(\tau)*\hat{v}_\omega^{\ast L-l}(\tau),\label{eqn:92}\\
\hat{m}_\omega^L(\tau)=\frac{1}{(2\pi)^{L-1}}\sum_{l=0}^L
\binom{L}{l}\left(\frac{e^{-(\omega+i\tau)}}{1-e^{-(\omega+i\tau)}}\right)^l
\left(\frac{e^{-\omega}}{1-e^{-\omega}}\right)^{L-l}
\hat{u}_\omega^{\ast l}(\tau)
*\hat{v}_\omega^{\ast L-l}(\tau),\label{eqn:93}
\end{align}
where $*$ denotes the convolution operation and $*l$ denotes the 
$l$-fold convolution.
\end{lemma}

\noindent{\bf Proof}\ \ See Appendix \ref{app:B}.\hfill$\Box$

\begin{lemma}
\label{lem:14}
For $\tau\neq0$, 
\begin{align}
\lim_{\omega\to0+}\hat{M}_\omega^L(\tau)&=
\frac{1}{(2\pi)^{L-1}}\sum_{l=0}^L\frac{1}{\pi^{L-l}}\binom{L}{l}
\left(\frac{1}{1-e^{-i\tau}}\right)^l\hat{q}_1^{*l}(\tau)
*\hat{q}_2^{*L-l}(\tau),\label{eqn:94}\\
\lim_{\omega\to0+}\hat{m}_\omega^L(\tau)&=
\frac{1}{(2\pi)^{L-1}}\sum_{l=0}^L\frac{1}{\pi^{L-l}}\binom{L}{l}
\left(\frac{e^{-i\tau}}{1-e^{-i\tau}}\right)^l\hat{q}_1^{*l}(\tau)
*\hat{q}_2^{*L-l}(\tau).\label{eqn:95}
\end{align}
\end{lemma}

\noindent{\bf Proof}\ \ (\ref{eqn:94}) and (\ref{eqn:95}) follow 
\begin{align}
\hat{u}_\omega^{*l}(\tau)&=\sum_{j=0}^l\left(-\frac{\omega}{\pi}\right)^{l-j}
\hat{q}_1^{*j}(\tau)*\hat{q}_2^{*l-j}(\tau),\label{eqn:96}\\
\hat{v}_\omega^{*L-l}(\tau)&=\left(\frac{\omega}{\pi}\right)^{L-l}
\hat{q}_2^{*L-l}(\tau).\label{eqn:97}
\end{align}

By the change of variables, we have

\begin{lemma}
\label{lem:15}
For $\tau\neq0$, 
\begin{align}
\lim_{\sigma\to0+}\hat{M}_{\sigma,\delta}^L(\tau)&=
\frac{2\pi}{\delta}\frac{1}{(2\pi)^{L-1}}\sum_{l=0}^L\frac{1}{\pi^{L-l}}\binom{L}{l}
\left(\frac{1}{1-e^{-i2\pi\tau/\delta}}\right)^l\hat{q}_1^{*l}\left(\frac{2\pi\tau}{\delta}\right)
*\hat{q}_2^{*L-l}\left(\frac{2\pi\tau}{\delta}\right),\label{eqn:98}\\
\lim_{\sigma\to0+}\hat{m}_{\sigma,\delta}^L(\tau)&=
\frac{2\pi}{\delta}\frac{1}{(2\pi)^{L-1}}\sum_{l=0}^L\frac{1}{\pi^{L-l}}\binom{L}{l}
\left(\frac{e^{-i2\pi\tau/\delta}}{1-e^{-i2\pi\tau/\delta}}\right)^l\hat{q}_1^{*l}
\left(\frac{2\pi\tau}{\delta}\right)
*\hat{q}_2^{*L-l}\left(\frac{2\pi\tau}{\delta}\right).\label{eqn:99}
\end{align}
\end{lemma}

\section{Proof of Theorem \ref{theo:1}}
\subsection{Upper Bound for $P(X>x)$}
First, we will evaluate $P(X>x)$ from above by using the 
majorant function $M_{\sigma,\delta}^L$. 

Let $L\in{\mathbb N}^+$ with $L\geq r$. For arbitrary 
$\sigma_1>0$, $\sigma_2>0$, $\delta>0$, 
\begin{align}
e^{L\sigma_2x}\int_x^\infty e^{-(\sigma_1+L\sigma_2)t}dF(t)
&=\int_x^\infty e^{-L\sigma_2(t-x)}e^{-\sigma_1t}dF(t)\label{eqn:100}\\
&=\int_0^\infty E_{L\sigma_2}(t-x)e^{-\sigma_1t}dF(t)\label{eqn:101}\\
&\leq\int_0^\infty M_{\sigma_2,\delta}^L(t-x)e^{-\sigma_1t}dF(t),\ x>0,\label{eqn:102}
\end{align}
where the last inequality holds by (\ref{eqn:79}) in Lemma \ref{lem:10} 
(see also Korevaar\,\cite{kor}, Nakagawa\,\cite{nak4}). 
By Lemma \ref{lem:11}, $M_{\sigma_2,\delta}^L(t-x)$ is represented by 
the inverse Fourier transform of 
$\hat{M}_{\sigma_2,\delta}^L={\cal F}\left(M_{\sigma_2,\delta}^L\right)$ as 
\begin{eqnarray}
M_{\sigma_2,\delta}^L(t-x)=\frac{1}{2\pi}\int_{-L\delta}^{L\delta}
\hat{M}_{\sigma_2,\delta}^L(-\tau)e^{i(x-t)\tau}d\tau.\label{eqn:103}
\end{eqnarray}
Substituting (\ref{eqn:103}) into (\ref{eqn:102}), we have 
by Fubini's theorem 
\begin{align}
\int_0^\infty M_{\sigma_2,\delta}^L(t-x)e^{-\sigma_1t}dF(t)
&=\int_0^\infty\frac{1}{2\pi}\int_{-L\delta}^{L\delta}\hat{M}_{\sigma_2,\delta}^L(-\tau)
e^{i(x-t)\tau}e^{-\sigma_1t}d\tau dF(t)\label{eqn:104}\\
&=\frac{1}{2\pi}\int_{-L\delta}^{L\delta}\hat{M}_{\sigma_2,\delta}^L(-\tau)
e^{ix\tau}d\tau\int_0^\infty e^{-(\sigma_1+i\tau)t}dF(t)\label{eqn:105}\\
&=\frac{1}{2\pi}\int_{-L\delta}^{L\delta}\hat{M}_{\sigma_2,\delta}^L(-\tau)
e^{ix\tau}\varphi(\sigma_1+i\tau)d\tau.\label{eqn:106}
\end{align}

Now, for the representation (\ref{eqn:16}) of $\varphi(s)$, 
i.e., 
\begin{eqnarray}
\varphi(s)=\alpha(s)s^r\log s +\beta(s),\ r\in{\mathbb N}^+,\label{eqn:107}
\end{eqnarray}
we define 
\begin{eqnarray}
f^\ast(t)=g^\ast(t)+h^\ast(t),\ t\geq0,\label{eqn:108}
\end{eqnarray}
where $g^\ast(t)$ and $h^\ast(t)$ were defined in Lemma \ref{lem:6}. 
Let $\varphi^\ast(s)={\cal L}(f^\ast(t))$, then $f^\ast(t)$ and 
$\varphi^\ast(s)$ have the following properties (a), (b), (c) and (c').

(a) $f^\ast(t)>0$ for all sufficiently large $t$.

\noindent This is because $g^\ast(0)=g_0>0$ by (\ref{eqn:58}).

Define $\xi(s)\equiv\varphi(s)-\varphi^\ast(s)$, then by Lemma 
\ref{lem:6}, 
\begin{align}
\xi(s)&=\varphi(s)-\left(G^\ast(s)+H^\ast(s)\right)\label{eqn:109}\\
&=s^L\left\{\left(\sum_{n=0}^\infty\alpha_{L+n}s^n\right)s^r\log s
+\sum_{n=0}^\infty\left(\beta_{L+n}-\tilde\beta_{L+n}\right)s^n\right\},\ 
s=\sigma+i\tau.\label{eqn:110}
\end{align}
We see from (\ref{eqn:110}),

(b) $\xi(s)$ is continuous in the closed region 
$\{0\leq\sigma\leq\epsilon,\ -L\delta\leq\tau\leq L\delta\}\subset{\mathbb C}$ 
for sufficiently small $\epsilon, \delta>0$.

\noindent In fact, because $s=0$ is an isolated singularity, 
we can take $\epsilon, \delta$ so small that the closed region 
$\{0\leq\sigma\leq\epsilon,\ -L\delta\leq\tau\leq L\delta\}$ 
does not include any singularities of $\varphi(s)$ and $\varphi^\ast(s)$ 
other than $s=0$.

Define 
$\hat{M}_{0,\delta}^L(\tau)=\lim_{\sigma\to0+}\hat{M}_{\sigma,\delta}^L(\tau),\ 
\hat{m}_{0,\delta}^L(\tau)=\lim_{\sigma\to0+}\hat{m}_{\sigma,\delta}^L(\tau),\,
-L\delta\leq\tau\leq L\delta,\,\tau\neq0$, for $\delta$ sufficiently small 
as in (b). Then, we have 

(c) $\mbox{\rm supp}(\hat{M}_{0,\delta}^L)\subset[-L\delta,\ L\delta]$ 
and $\hat{M}_{0,\delta}^L(\tau)\xi(i\tau)$ is $r$ times 
piecewise differentiable with 

\hspace*{5mm}$\left(\hat{M}_{0,\delta}^L(\tau)\xi(i\tau)\right)^{(r)}\in 
L^1\Big([-L\delta,\ L\delta]\Big)$.

\medskip

(c') $\mbox{\rm supp}(\hat{m}_{0,\delta}^L)\subset[-L\delta,\ L\delta]$ 
and $\hat{m}_{0,\delta}^L(\tau)\xi(i\tau)$ is $r$ times 
piecewise differentiable with 

\hspace*{5mm}$\left(\hat{m}_{0,\delta}^L(\tau)\xi(i\tau)\right)^{(r)}\in 
L^1\Big([-L\delta,\ L\delta]\Big)$.

\noindent These (c) and (c') hold because of Lemma \ref{lem:15} and 
(\ref{eqn:110}). 

Now, defining $F^\ast(t)$ by $dF^\ast(t)=f^\ast(t)dt$, we have 
in a similar way as (\ref{eqn:106}) 
\begin{eqnarray}
\int_0^\infty M_{\sigma_2,\delta}^L(t-x)e^{-\sigma_1t}dF^\ast(t)
=\frac{1}{2\pi}\int_{-L\delta}^{L\delta}\hat{M}_{\sigma_2,\delta}^L(-\tau)e^{ix\tau}
\varphi^\ast(\sigma_1+i\tau)d\tau.\label{eqn:111}
\end{eqnarray}
Subtracting (\ref{eqn:111}) from (\ref{eqn:106}), we have 
\begin{align}
\int_0^\infty M_{\sigma_2,\delta}^L(t-x)e^{-\sigma_1t}dF(t)
&=\int_0^\infty M_{\sigma_2,\delta}^L(t-x)e^{-\sigma_1t}dF^\ast(t)\label{eqn:112}\\
&\hspace{10mm}+\frac{1}{2\pi}\int_{-L\delta}^{L\delta}\hat{M}_{\sigma_2,\delta}^L(-\tau)e^{ix\tau}
\xi(\sigma_1+i\tau)d\tau.\label{eqn:113}
\end{align}
For sufficiently small $\delta>0$, 
\begin{eqnarray}
\xi(i\tau)=\lim_{\sigma_1\to0+}\xi(\sigma_1+i\tau),\ -L\delta\leq\tau\leq L\delta\label{eqn:114}
\end{eqnarray}
is uniform convergence due to (b), hence 
\begin{eqnarray}
\lim_{\sigma_1\to0+}\frac{1}{2\pi}\int_{-L\delta}^{L\delta}
\hat{M}_{\sigma_2,\delta}^L(-\tau)\xi(\sigma_1+i\tau)e^{ix\tau}d\tau
=\frac{1}{2\pi}\int_{-L\delta}^{L\delta}
\hat{M}_{\sigma_2,\delta}^L(-\tau)\xi(i\tau)e^{ix\tau}d\tau.\label{eqn:115}
\end{eqnarray}
From (\ref{eqn:102}), (\ref{eqn:113}), (\ref{eqn:115}), for 
$\sigma_1\to0+$, we have 
\begin{eqnarray}
e^{L\sigma_2 x}\int_x^\infty e^{-L\sigma_2 t}dF(t)
\leq\int_0^\infty M_{\sigma_2,\delta}^L(t-x)dF^\ast(t)
+\frac{1}{2\pi}\int_{-L\delta}^{L\delta}\hat{M}_{\sigma_2,\delta}^L(-\tau)
\xi(i\tau)e^{ix\tau}d\tau.\label{eqn:116}
\end{eqnarray}
By the estimation for $M_\omega^1(t)$ (see Korevaar\,\cite{kor}, p.132), 
i.e., 
\begin{eqnarray}
\left\{\begin{array}{ll}\label{eqn:117}
0\leq M_\omega^1(t)\leq\left(\ds\frac{\sin\pi t}{\pi t}\right)^2, & t<0,\\
e^{-\omega t}\leq M_\omega^1(t)\leq e^{-\omega t}+\left(\ds\frac{\sin\pi t}{\pi t}\right)^2, 
& t\geq0,
\end{array}\right.
\end{eqnarray}
we have
\begin{eqnarray}
\left\{\begin{array}{ll}\label{eqn:118}
0\leq M_{\sigma,\delta}^L(t)\leq\left(\ds\frac{\sin\delta t/2}{\delta t/2}\right)^{2L}, & t<0,\\
e^{-L\omega t}\leq M_{\sigma,\delta}^L(t)\leq\left(e^{-\omega t}+\left(\ds\frac{\sin\pi t}{\pi t}\right)^2\right)^L, & t\geq0,
\end{array}\right.
\end{eqnarray}
where $\omega=2\pi\sigma/\delta$. Thus, it is easy to see that 
there exists a constant $C_1>0$ such that
\begin{eqnarray}
M_{\sigma,\delta}^L(t-x)\leq\left\{\begin{array}{ll}\label{eqn:119}
\left(\ds\frac{1}{\delta(x-t)/2}\right)^{2L}, & \ 0\leq t<x-1,\\[5mm]
C_1, & t\geq x-1.
\end{array}\right.
\end{eqnarray}
Therefore, the first term of the right hand side of (\ref{eqn:116}) 
is evaluated as 
\begin{align}
\int_0^\infty M_{\sigma_2,\delta}^L&(t-x)dF^\ast(t)
=\int_0^\infty M_{\sigma_2,\delta}^L(t-x)f^\ast(t)dt\label{eqn:120}\\
&=\int_0^\infty M_{\sigma_2,\delta}^L(t-x)\left(g^\ast(t)+h^\ast(t)\right)dt\label{eqn:121}\\
&\leq\sum_{k=0}^{L-1}|g_k|
\int_1^{x-1}\left(\ds\frac{1}{\delta(x-t)/2}\right)^{2L}
\frac{1}{t^{r+k+1}}dt+C_1\sum_{k=0}^{L-1}|g_k|\int_{x-1}^\infty
\ds\frac{1}{t^{r+k+1}}dt\label{eqn:122}\\
&\hspace{5mm}+\sum_{k=1}^Lk|d_k|
\int_0^{x-1}\left(\ds\frac{1}{\delta(x-t)/2}\right)^{2L}e^{-kt}dt+C_1\sum_{k=1}^Lk|d_k|\int_{x-1}^\infty e^{-kt}dt\label{eqn:123}\\
&\leq O(x^{-(r+1)})+C_2g_0x^{-r}+O(x^{-2L})+O(e^{-x}),\ C_2>0\label{eqn:124}\\
&<\ds\frac{C_3}{x^r},\ C_3>0,\label{eqn:125}
\end{align}
for all sufficiently large $x$, by virtue of Lemmas \ref{lem:A1}, 
\ref{lem:A2} in Appendix. Notice $g_0>0$ and $L\geq r$.

Next, the second term of the right hand side of (\ref{eqn:116}) 
will be estimated. We have by (c) and integration by parts, 
\begin{align}
\lim_{\sigma_2\to0+}\ds\frac{1}{2\pi}\int_{-L\delta}^{L\delta}\hat{M}_{\sigma_2,\delta}^L
(-\tau)&\xi(i\tau)e^{ix\tau}d\tau
=\ds\frac{1}{2\pi}\int_{-L\delta}^{L\delta}\hat{M}_{0,\delta}^L
(-\tau)\xi(i\tau)e^{ix\tau}d\tau\label{eqn:126}\\
&=\ds\frac{i^r}{2\pi x^r}\int_{-L\delta}^{L\delta}\left(\hat{M}_{0,\delta}^L
(-\tau)\xi(i\tau)\right)^{(r)}e^{ix\tau}d\tau\label{eqn:127}\\
&=o(x^{-r}),\ x\to\infty,\label{eqn:128}
\end{align}
due to Riemann-Lebesgue theorem. Then in (\ref{eqn:116}) 
for $\sigma_2\to0+$, we have by (\ref{eqn:125}) and (\ref{eqn:128}), 
\begin{eqnarray}
P(X>x)=\int_x^\infty dF(t)<\ds\frac{C}{x^r},\ C>0,\label{eqn:129}
\end{eqnarray}
for all sufficiently large $x$.

\subsection{Lower Bound for $P(X>x)$}
We will evaluate $P(X>x)$ from below by using the minorant function $m_{\sigma,\delta}^L$.

Let $L\in{\mathbb N}^+$ be an odd number with $L\geq r$. 
For arbitrary $\sigma_1>0$, $\sigma_2>0$, $\delta>0$, 
\begin{align}
e^{L\sigma_2 x}\int_x^\infty e^{-(\sigma_1+L\sigma_2)t}dF(t)
&=\int_0^\infty E_{L\sigma_2}(t-x)e^{-\sigma_1t}dF(t)\label{eqn:130}\\
&\geq\int_0^\infty m_{\sigma_2,\delta}^L(t-x)e^{-\sigma_1 t}dF(t),\ x>0.\label{eqn:131}
\end{align}
In a similar way as from (\ref{eqn:100}) to (\ref{eqn:116}), 
we have
\begin{eqnarray}
e^{L\sigma_2 x}\int_x^\infty e^{-L\sigma_2 t}dF(t)
\geq\int_0^\infty m_{\sigma_2,\delta}^L(t-x)dF^\ast(t)
+\ds\frac{1}{2\pi}\int_{-L\delta}^{L\delta}\hat{m}_{\sigma_2,\delta}^L
(-\tau)\xi(i\tau)e^{ix\tau}d\tau.\label{eqn:132}
\end{eqnarray}
Due to the estimation for $m_\omega^1(t)$ (see Korevaar\,\cite{kor}, p.132), 
i.e., 
\begin{eqnarray}
\left\{\begin{array}{ll}\label{eqn:133}
-\left(\ds\frac{\sin\pi t}{\pi t}\right)^2\leq m_\omega^1(t)\leq0, & t<0,\\
e^{-\omega t}-\left(\ds\frac{\sin\pi t}{\pi t}\right)^2\leq m_\omega^1(t)
\leq e^{-\omega t} & t\geq0,
\end{array}\right.
\end{eqnarray}
we have
\begin{eqnarray}
\left\{\begin{array}{ll}\label{eqn:134}
-\left(\ds\frac{\sin\delta t/2}{\delta t/2}\right)^{2L}\leq m_{\sigma,\delta}^L(t)
\leq0, & t<0,\\
\left(e^{-\omega t}-\left(\ds\frac{\sin\delta t/2}{\delta t/2}\right)^2\right)^L
\leq m_{\sigma,\delta}^L(t)\leq e^{-L\omega t}, & t\geq0,
\end{array}\right.
\end{eqnarray}
where $\omega=2\pi\sigma/\delta$. Thus, there exit 
constants $C_4, C_5, C_6>0$ such that
\begin{eqnarray}
m_{\sigma,\delta}^L(t-x)\geq\left\{\begin{array}{ll}\label{eqn:135}
-\left(\ds\frac{1}{\delta(x-t)/2}\right)^{2L}, & \ 0\leq t<x-1,\\[5mm]
-C_4, & x-1\leq t<x+1,\\[3mm]
C_5e^{-(2\pi L\sigma/\delta)(t-x)}-C_6\left(\ds\frac{1}{\delta(x-t)/2}\right)^2,
& t\geq x+1.
\end{array}\right.
\end{eqnarray}
Therefore, the first term of the right hand side of (\ref{eqn:132}) 
is evaluated as 
\begin{align}
\int_0^\infty m_{\sigma_2,\delta}^L(t-x)(g^\ast(t)&+h^\ast(t))dt\geq
-\sum_{k=0}^{L-1}|g_k|\int_1^{x-1}\left(\frac{1}{\delta(x-t)/2}\right)^{2L}
\frac{1}{t^{r+k+1}}dt\label{eqn:136}\\
&-C_4\sum_{k=0}^{L-1}|g_k|\int_{x-1}^{x+1}\frac{1}{t^{r+k+1}}dt\label{eqn:137}\\
&+C_5\sum_{k=0}^{L-1}|g_k|\int_{x+1}^\infty e^{-(2\pi L\sigma_2/\delta)(t-x)}
\frac{1}{t^{r+k+1}}dt\label{eqn:138}\\
&-C_6\sum_{k=0}^{L-1}|g_k|\int_{x+1}^\infty \left(\frac{1}{\delta(x-t)/2}\right)^2
\frac{1}{t^{r+k+1}}dt\label{eqn:139}\\
&-\sum_{k=1}^Lk|d_k|\int_0^{x-1}\left(\frac{1}{\delta(x-t)/2}\right)^{2L}
e^{-kt}dt\label{eqn:140}\\
&-C_4\sum_{k=1}^Lk|d_k|\int_{x-1}^{x+1}e^{-kt}dt\label{eqn:141}\\
&+C_5\sum_{k=1}^Lk|d_k|\int_{x+1}^\infty e^{-(2\pi L\sigma/\delta)(t-x)}
e^{-kt}dt\label{eqn:142}\\
&-C_6\sum_{k=1}^Lk|d_k|\int_{x+1}^\infty\left(\frac{1}{\delta(x-t)/2}\right)^2
e^{-kt}dt\label{eqn:143}\\
&\geq O(x^{-(r+1)})+O(x^{-(r+1)})+C_5'g_0\int_{x+1}^\infty e^{-(2\pi L\sigma_2/\delta)(t-x)}
\frac{1}{t^{r+1}}dt\label{eqn:144}\\
&\hspace{5mm}+O(x^{-(r+1)})+O(x^{-2L})+O(e^{-x})+O(e^{-x})+O(e^{-x}),\label{eqn:145}
\end{align}
where $C_5'$ is a positive constant. Thus, 
\begin{eqnarray}
\lim_{\sigma_2\to0+}\int_0^\infty m_{\sigma_2,\delta}^L(t-x)dF^\ast(t)
\geq\ds\frac{C_7}{x^r},\ C_7>0,\label{eqn:146}
\end{eqnarray}
for all sufficiently large $x$.

Next, we will evaluate the second term of the right hand 
side of (\ref{eqn:132}). We have, by (c') and the integration 
by parts, 
\begin{align}
\lim_{\sigma_2\to0+}\ds\frac{1}{2\pi}\int_{-L\delta}^{L\delta}
\hat{m}_{\sigma_2,\delta}^L(-\tau)\xi(i\tau)e^{ix\tau}d\tau
&=\frac{1}{2\pi}\int_{-L\delta}^{L\delta}
\hat{m}_{0,\delta}^L(-\tau)\xi(i\tau)e^{ix\tau}d\tau\label{eqn:147}\\
&=\frac{i^r}{2\pi x^r}\int_{-L\delta}^{L\delta}
\Big(\hat{m}_{0,\delta}^L(-\tau)\xi(i\tau)\Big)^{(r)}e^{ix\tau}d\tau\label{eqn:148}\\
&=o(x^{-r}),\ x\to\infty,\label{eqn:149}
\end{align}
by Riemann-Lebesgue theorem. Then, in (\ref{eqn:132}) for 
$\sigma_2\to0+$, we have from (\ref{eqn:146}) and (\ref{eqn:149}), 
\begin{eqnarray}
P(X>x)=\int_x^\infty dF(t)>\frac{C'}{x^r},\ C'>0,\label{eqn:150}
\end{eqnarray}
for all sufficiently large $x$.

From (\ref{eqn:129}) and (\ref{eqn:150}), the proof of 
Theorem \ref{theo:1} is completed. \hfill$\Box$

\section{Proof of Theorem \ref{theo:2}}
The same proof as that of Theorem \ref{theo:1} is applicable by 
replacing $g^\ast(t)$, $h^\ast(t)$ in Lemma \ref{lem:5} 
with those in Lemma \ref{lem:6}.\hfill$\Box$

\section{Conclusion}
In this paper, we investigated the asymptotic decay of the 
tail probability of a heavy tailed random variable. 
We proved two theorems which give sufficient conditions 
for a random variable to be heavy tailed.
Our theorems are 
based on the Tauberian theorems due to Graham and Vaarler. 
The central idea is the approximation of the exponential 
function by a majorant and minorant functions whose Fourier 
transforms have a finite support. 

Through the proof of the theorems, I think that some more 
general representation of the singularity guarantees the 
random variable to be heavy tailed.


\newpage


\newpage

\appendix

\noindent{\Large\bf Appendix}

\begin{lemma}
\label{lem:A1}
For $n_1, n_2\in{\mathbb N}^+,\ n_1\geq2,\ n_2\geq2$, let 
$n=\min(n_1,n_2)$. Then,
\begin{eqnarray*}
\int_1^{x-1}\ds\frac{dt}{(x-t)^{n_1}t^{n_2}}\leq O(x^{-n}),\ x\to\infty.
\end{eqnarray*}
\end{lemma}

\noindent{\bf Proof}\ \ By the change of variable $t=xu$,
\begin{align*}
\int_1^{x-1}\ds\frac{dt}{(x-t)^{n_1}t^{n_2}}
&=\ds\frac{1}{x^{n_1+n_2+1}}
\int_{1/x}^{1-1/x}\ds\frac{du}{(1-u)^{n_1}u^{n_2}}\\
&=\ds\frac{1}{x^{n_1+n_2+1}}\left\{\int_{1/x}^{1/2}+\int_{1/2}^{1-1/x}\right\}
\frac{du}{(1-u)^{n_1}u^{n_2}}\\
&\leq\ds\frac{1}{x^{n_1+n_2+1}}\left\{2^{n_1}\int_{1/x}^{1/2}\frac{du}{u^{n_2}}
+2^{n_2}\int_{1/2}^{1-1/x}\frac{du}{(1-u)^{n_1}}\right\}\\
&=\ds\frac{1}{x^{n_1+n_2+1}}\left\{2^{n_1}\frac{x^{n_2-1}-2^{n_2-1}}{n_2-1}
+2^{n_2}\frac{x^{n_1-1}-2^{n_1-1}}{n_1-1}\right\}\\
&=O(x^{-n}),\ x\to\infty.
\end{align*}

\begin{lemma}
\label{lem:A2}
For $k>0$ and $n\in{\mathbb N}$, we have
\begin{eqnarray*}
\int_0^{x-1}\ds\frac{e^{-kt}}{(x-t)^n}dt=O(x^{-n}),\ x\to\infty.
\end{eqnarray*}
\end{lemma}

\noindent{\bf Proof}\ \ By the change of variable $u=x-t$, 
\begin{align*}
\int_0^{x-1}\ds\frac{e^{-kt}}{(x-t)^n}dt&=e^{-kx}\int_1^x\ds\frac{e^ku}{u^n}du\\
&=\ds\frac{1}{x^n}\int_1^x\frac{e^{ku}}{u^n}du\Big/\frac{e^{kx}}{x^n}\\
&\rightarrow\ds\frac{1}{kx^n},\ x\to\infty,
\end{align*}
from L'Hopital's rule. \hfill$\Box$

\section{Proof of Lemma \ref{lem:12}}
\label{app:A}
\begin{align*}
\hat{M}_\omega^1(\tau)&=\sum_{n=0}^\infty e^{-n\omega}
{\cal F}\left(\left(\ds\frac{\sin\pi(t-n)}{\pi(t-n)}\right)^2\right)
-\frac{\omega}{\pi}\sum_{n=0}^\infty e^{-n\omega}
\left\{{\cal F}\left(\frac{\sin^2\pi(t-n)}{\pi(t-n)}\right)
-{\cal F}\left(\frac{\sin^2\pi t}{\pi t}\right)\right\}\\
&=\hat{q}_1(\tau)\sum_{n=0}^\infty e^{-n(\omega+i\tau)}
-\frac{\omega}{\pi}\hat{q}_2(\tau)\left(\sum_{n=0}^\infty e^{-n(\omega+i\tau)}
-\sum_{n=0}^\infty e^{-n\omega}\right)\\
&=\frac{1}{1-e^{-(\omega+i\tau)}}\hat{q}_1(\tau)
-\frac{\omega}{\pi}\left(\frac{1}{1-e^{-(\omega+i\tau)}}
-\frac{1}{1-e^{-\omega}}\right)\hat{q}_2(\tau).
\end{align*}
The result for $\hat{m}_\omega^1$ is proved in a similar way. \hfill$\Box$

\section{Proof of Lemma \ref{lem:13}}
\label{app:B}
\begin{align*}
M_\omega^L(t)&=\left(\frac{\sin\pi t}{\pi t}\right)^{2L}\left(Q_\omega(t)\right)^L\\
&=\left(\frac{\sin\pi t}{\pi t}\right)^{2L}\sum_{n_1=0}^\infty\ldots\sum_{n_L=0}^\infty 
e^{-n_1\omega}\ldots e^{-n_L\omega}\left\{\frac{1}{(t-n_1)^2}-\frac{\omega}{t-n_1}
+\frac{\omega}{t}\right\}\times\ldots\\
&\hspace{30mm}\times\left\{\frac{1}{(t-n_L)^2}-\frac{\omega}{t-n_L}
+\frac{\omega}{t}\right\}\\
&=\sum_{n_1=0}^\infty\ldots\sum_{n_L=0}^\infty e^{-n_1\omega}\ldots e^{-n_L\omega}
\left\{\left(\frac{\sin\pi(t-n_1)}{\pi(t-n_1)}\right)^2-\frac{\omega}{\pi}\frac{\sin^2\pi(t-n_1)}{\pi(t-n_1)}+\frac{\omega}{\pi}\frac{\sin^2\pi t}{\pi t}\right\}\times\ldots\\
&\hspace{30mm}
\times\left\{\left(\frac{\sin\pi(t-n_L)}{\pi(t-n_L)}\right)^2-\frac{\omega}{\pi}\frac{\sin^2\pi(t-n_L)}{\pi(t-n_L)}+\frac{\omega}{\pi}\frac{\sin^2\pi t}{\pi t}\right\}\\
&=\sum_{n_1=0}^\infty\ldots\sum_{n_L=0}^\infty e^{-n_1\omega}\ldots e^{-n_L\omega}
\left\{u_\omega(t-n_1)+v_\omega(t)\right\}\times\ldots
\times\left\{u_\omega(t-n_L)+v_\omega(t)\right\}.
\end{align*}
Therefore,
\begin{align}
\hat{M}_\omega^L(\tau)&=\sum_{n_1=0}^\infty\ldots\sum_{n_L=0}^\infty e^{-n_1\omega}
\ldots e^{-n_L\omega}
\frac{1}{(2\pi)^{L-1}}\left\{e^{-in_1\tau}\hat{u}_\omega(\tau)+\hat{v}_\omega(\tau)\right\}
*\ldots*\left\{e^{-in_L\tau}\hat{u}_\omega(\tau)+\hat{v}_\omega(\tau)\right\}\nonumber\\
&=\frac{1}{(2\pi)^{L-1}}\sum_{n_1=0}^\infty\ldots\sum_{n_L=0}^\infty e^{-n_1\omega}
\ldots e^{-n_L\omega}\nonumber\\
&\hspace{10mm}\times\sum_{l=0}^L\left\{\sum_{k_1,\ldots,k_l\mbox{\rm\,:\,distinct}}
\left(e^{-in_{k_1}\tau}\hat{u}_\omega(\tau)\right)*\ldots*
\left(e^{-in_{k_l}\tau}\hat{u}_\omega(\tau)\right)\right\}*\hat{v}_\omega^{*L-l}(\tau)
\label{eqn:distinct1}\\
&=\frac{1}{(2\pi)^{L-1}}\sum_{n_1=0}^\infty\ldots\sum_{n_L=0}^\infty e^{-n_1\omega}
\ldots e^{-n_L\omega}\nonumber\\
&\hspace{10mm}\times\sum_{l=0}^L\left(\sum_{k_1,\ldots,k_l\mbox{\rm\,:\,distinct}}
e^{-in_{k_1}\tau}\ldots e^{-in_{k_l}\tau}\right)
\hat{u}_\omega^{*l}(\tau)*\hat{v}_\omega^{*L-l}(\tau)
\label{eqn:distinct2}\\
&=\frac{1}{(2\pi)^{L-1}}\sum_{l=0}^L\hat{u}_\omega^{*l}(\tau)*\hat{v}_\omega^{*L-l}(\tau)\nonumber\\
&\hspace{10mm}
\times\sum_{n_1=0}^\infty\ldots\sum_{n_L=0}^\infty e^{-n_1\omega}\ldots e^{-n_L\omega}
\left(\sum_{k_1,\ldots,k_l\mbox{\rm\,:\,distinct}}
e^{-in_{k_1}\tau}\ldots e^{-in_{k_l}\tau}\right)
\label{eqn:distinct3}\\
&=\frac{1}{(2\pi)^{L-1}}\sum_{l=0}^L\hat{u}_\omega^{*l}(\tau)*\hat{v}_\omega^{*L-l}(\tau)
\binom{L}{l}\left(\frac{1}{1-e^{-(\omega+i\tau)}}\right)^l
\left(\frac{1}{1-e^{-\omega}}\right)^{L-l}\nonumber\\
&=\frac{1}{(2\pi)^{L-1}}\sum_{l=0}^L\binom{L}{l}\left(\frac{1}{1-e^{-(\omega+i\tau)}}\right)^l
\left(\frac{1}{1-e^{-\omega}}\right)^{L-l}
\hat{u}_\omega^{*l}(\tau)*\hat{v}_\omega^{*L-l}(\tau).\nonumber
\end{align}
In (\ref{eqn:distinct1}),(\ref{eqn:distinct2}),(\ref{eqn:distinct3}), for 
$l=0$, the (empty) sum on $k_1,\cdots,k_l$ is considered to be 1.

Similarly, we have the result for $\hat{m}_\omega^L(\tau)$.
\end{document}